\documentclass[10pt,amssymb,amscd]{amsart}
\usepackage{color}
\usepackage{amssymb}
\usepackage{amsmath}
\usepackage{amsthm} 
\usepackage{mathrsfs}                                                           
\usepackage{slashbox}  
\usepackage[all]{xy}
\usepackage{graphicx}
\usepackage{wrapfig}

\def\2{{1\over 2}}

\def\bb{{\bf b}}

\newcommand{\rf}[1]{(\ref{#1})}
\def\b{\bar}

\newcommand{\ud}{\mathrm{d}}

\def\b{\bar}

\def\<{\langle}
\def\>{\rangle}
\def\+{\dagger}

\begin{document}
\title{On higher order Leibniz identities in TCFT}
\author{Anton M. Zeitlin}
\address{ 
\newline Department of Mathematics,\newline
Columbia University,
\newline 2990 Broadway, New York,\newline
NY 10027, USA.
\newline
zeitlin@math.columbia.edu\newline
http://math.columbia.edu/$\sim$zeitlin \newline
http://www.ipme.ru/zam.html  }

\begin{abstract}
We extend the algebra of local observables in topological conformal field theories by nonlocal operators. This allows to construct parameter-dependent operations realized via certain integrals over the compactified moduli spaces, satisfying analogues of the Leibniz and higher order Leibniz identities holding up to homotopy. We conjecture that one can construct a complete set of such operations which lead to a parameter-dependent version of Loday's homotopy Leibniz algebras. 
\end{abstract}
\maketitle
\section{Introduction}
Understanding of the underlying hidden symmetries of two-dimensional topological conformal field theory (TCFT) is an important problem, both because it is a useful playground for higher-dimensional topological physical theories and because TCFT gives the proper formulation of the String theory.   
Usually it is useful to decompose TCFT into chiral and antichiral sectors and study both of them separately. 
The simplest examples of such (anti)chiral sectors are described by a mathematical object known as vertex operator algebra (VOA) \cite{benzvi} which in the case of topological CFT is called topological VOA (TVOA).  

Striking relation between TVOAs and homotopy Gerstenhaber ($G_{\infty}$) homotopy algebras was proposed in \cite{lz} and then proved in 
 \cite{zuck},\cite{huang}. That relation was also extended in \cite{gorbounov}, \cite{vallette} to the case of $BV_{\infty}$ algebras. 

The objects of interest in this article are the symmetries which occur in the case of full TCFT. In \cite{zeitlin} we  considered operators in TCFT on the real line only and extended the space of states by nonlocal operators, which were obtained by means of integration of the correlation function over compact manifolds with boundary, embedded in a simplex in $\mathbb{R}^n$. We conjectured that there is a structure of parameter-dependent $A_{\infty}$ algebra \cite{stasheff},\cite{stashbook} on the resulting new space of nonlocal operators, where the operations are obtained by means of integration over Stasheff polytopes. The $A_{\infty}$ algebra relations were checked up to "pentagon" equation involving bilinear and trilinear ones. We also noted the relation of this algebra to the compactified real slices of moduli space of points on the real line \cite{kap}, \cite{devadoss}, based on heuristic arguments of physics papers (see e.g. \cite{dijk}, \cite{herbst}, \cite{khromov} and references therein).

In this paper we are working with full TCFT on the whole complex plane and we extend the space of states by means of nonlocal operators, which are constructed by means of the integration of correlation functions over a compact manifold with boundary in $\mathbb{C}^n$. It is possible to extend this space by tensoring it with the 2d de Rham complex and as a result  we obtain the space of nonlocal operator-valued 0-, 1- and 2-forms. We claim that the resulting space has a system of n-linear operations on it, satisfying a  parameter-dependent version of Loday's homotopy Leibniz algebra $(Leibniz_{\infty})$ \cite{loday}, \cite{xu} up to the homotopy with respect to de Rham differential. We construct these 
operations up to $n=4$ and we notice that they are described by means of integration over compactified moduli spaces $\bar{M}(n)$ (see e.g. \cite{zuck},\cite{konts}, \cite{merkulov} and references therein). Based on that, we conjecture the explicit formula for higher operations.

This algebra is very similar to $L_{\infty}$ algebra which makes possible the relation between this algebra and $L_{\infty}$ algebra of String Field Theory \cite{zwiebach}, where the construction of higher order operations is somewhat similar. 

The structure of the paper is as follows. In Section 2 we describe the Leibniz differential algebra, related to TVOA, where we extend the space of TVOA by 1-forms. The restriction of this algebra to 0-forms gives the Leibniz algebra of Lian and Zuckerman \cite{lz}. 

In Section 3 we define a space of nonlocal operators and then construct first four operations satisfying $Leibniz_{\infty}$ relations up to homotopy with respect to de Rham operator. 

In Section 4 we present the formula for higher n-linear operations motivated by the results of Section 3. In the end we discuss the physical implications of our construction, in particular to the study of conformal perturbation theory.

\section{Toy story: Leibniz algebras for TVOA}
\noindent{\bf 2.0. Notation.} We assume that the reader is familiar with vertex operator algebras (VOA). For simplicity, we denote the vertex operator associated to the vector $A\in V$, where $V$ is the VOA space, as $A(z)$ instead of usual designation $Y(A,z)$. We will often refer to vectors in the space $V$ of the VOA as $states$.\\

\noindent {\bf 2.1. Leibniz algebra.} 
In this section we will describe Leibniz algebra associated with any topological vertex operator algebra (TVOA). In section 3 this construction will be generalized to the case of full Topological CFT. 

First we need to recall the definition of TVOA.\\

\noindent{\bf Definition 2.1.} {\it Let V be a $\mathbb{Z}$-graded vertex operator superalgebra, such that $V=\oplus_{i\in \mathbb{Z}} V^i=\oplus_{i,\mu\in \mathbb{Z}}V^i[\mu]$, where $i$ represents grading of $V$ with respect to conformal weight and $\mu$ 
represents fermionic grading of $V^i$. We call V a topological vertex operator algebra (TVOA) if there exist four elements: $J\in V^1[1]$, $b\in V^2[-1]$, $N\in V^1[0]$, $L\in V^2[0]$, such that 
\begin{eqnarray}
[Q,b(z)]=\mathcal{L}(z),\quad  Q^2=0,\quad b_0^2=0,
\end{eqnarray}
where $Q=J_0$ and $b(z)=\sum_n b_nz^{-n-2}$, $J(z)=\sum_n J_nz^{-n-1}$, \\
$\mathcal{L}(z)=\sum_n\mathcal{L}_nz^{-n-2}$, 
$F(z)=\sum_nF_nz^{-n-1}$. 
Here $\mathcal{L}(z)$ is the Virasoro element of $V$; the operators $N_0$, $\mathcal{L}_0$ are diagonalizable, commute with each other and their eigenvalues 
coincide with fermionic grading and conformal weight correspondingly.}\\

Let us introduce the operator $s$ which shifts the degree by +1. Then one can consider the space $W=V\oplus V'$, so that $V'\cong sV$. There are two maps $\ud: W\to W$ of degree 1 and $b: W\to W$ of degree 0. 
The operator $\ud:V\to V'$  is given by $sL_{-1}(-1)^{N_0}$ and zero when applied to $V'$, while  $b:V\to V'$ is given by $s(-1)^{N_0}b_{-1}$ and zero when applied to $V'$. 

It is more natural to interpret these operators if we consider the language of vertex operators. For a formal varIable $z$, let us introduce an odd formal variable $\theta$ of degree 1, so that $\theta^2=0$. 
If a vertex operator $A(z)$ corresponds to the element $A\in V$, let us associate to $sA\in V'$ the following object, which we will call a 1-form: 
$A(z)\theta$, so that $sA=(-1)^{|A|}s\frac{d}{d\theta}\lim_{z\to 0}A(z)\theta|0\rangle$. Obviously, $A(z)\theta=(-1)^{|A|}\theta A(z)$.  According to this procedure, one can associate to any element $\alpha \in W$ of a given degree a formal power series $\alpha(z,\theta)\in End(W)[[z, z^{-1}]][\theta]$, so that $\alpha=\lim_{z, \theta\to 0}(1-(-1)^{|\alpha|} s\frac{d}{d\theta}\alpha(z,\theta)|0\rangle$. 

On the level of vertex operators and their 1-forms, the operator $\ud$ becomes a de Rham operator $\ud=\theta L_{-1}$ and the operator $b$ acts as follows: $\alpha(z)\to [\alpha(z), {\bf b}]$, where 
${\bf b}(\theta)=b_{-1}\theta$. 

As usual in the theory of vertex algebras, 
under the correlator formal variable $z$ becomes a complex variable. Therefore, under the correlator one can integrate 1-forms we introduced, i.e.
$\int_{C}\phi=\int_{C}A(z)dz$, if $\phi=A(z)\theta$.

This allows us to define the following bilinear operation on the space $W$, which looks as follows when considered under the correlator:
\begin{eqnarray}\label{circle}
\alpha\circ \beta(z,\theta)=\int_{C_z} \alpha(w, \theta')\beta(z,\theta)+\int_{C_z}[\alpha(w,\theta'),{\bf b}(\theta')]\beta(z,\theta),
\end{eqnarray}
where $\alpha, \beta\in W$ the integral is over $w$-variable and $C_z$ is a closed contour around $z$. One can notice that if $\alpha\in V$ then only the second term of \rf{circle} is nonzero and if $\alpha\in V'$ only the first term survives. 
Using this bilinear operation it is possible to define another bilinear operation of degree $-1$ (which as we show below, satisfies the Leibniz algebra relations):
\begin{eqnarray}\label{bracket}
[\alpha,\beta](z,\theta)=\alpha\circ\beta(z, \theta)-[\alpha, {\bf b}(\theta)]\circ\beta(z,\theta).
\end{eqnarray}
Again, we note that because of $\theta^2=0$, the second term contributes only when $\beta\in V$.  
Let us return back to the definition of TVOA and notice that due to the presence of the nilpotent operator $Q$ of degree +1, which acts on both $V$ and $V'$, $W$ is actually a bicomplex. The differentials are related via familiar formula: $[Q, {\bf b}]=\ud$. 
Let us introduce another differential which will be relevant in the following:
\begin{eqnarray}
\bar{D}=Q-\ud.
\end{eqnarray}
Now we are ready to formulate a proposition.\\
 
 \noindent{\bf Theorem 2.1.}{\it  The bilinear operation $[\cdot, \cdot]$ of degree $-1$ satisfies the relations of  differential Leibniz algebra on $W$:
\begin{eqnarray}\label{leib}
&&\bar{D}[\alpha, \beta]=[\bar{D}\alpha, \beta]+(-1)^{|\alpha|+1}[\alpha, \bar{D}\beta],\nonumber\\
&&[\alpha, [\beta, \gamma]]=[[\alpha, \beta], \gamma]+(-1)^{(|\alpha|+1)(|\beta|+1)}[\beta, [\alpha, \gamma]].
\end{eqnarray}
  }\\
  
\noindent{\bf Proof.} 
Let us prove the first relation. 
\begin{eqnarray}
&&\bar{D}[\alpha, \beta]=(Q-\ud)[\alpha, \beta]=\nonumber\\
&&Q\alpha\circ \beta+(-1)^{|\alpha|+1}\alpha\circ Q\beta-\ud(\alpha\circ\beta)+[\ud,\alpha]\circ \beta - \nonumber\\
&&[Q\alpha,{\bf b}]\circ \beta-
(-1)^{|\alpha|+1}[\alpha,{\bf b}]\circ Q\beta=\nonumber \\
&&[\bar{D}\alpha, \beta]+(-1)^{|\alpha|+1}[\alpha, \bar{D}\beta].
\end{eqnarray}
Here we have applied shorthand notation,  $[A, {\bf b}]\circ B$ (resp.  $[A, {\bf d}]\circ B$) 
for the vector corresponding to the operator  $[A, {\bf b}(\theta)]\circ B(z,\theta)$ (resp.  $[A, \ud_{\theta}]\circ B(z,\theta)$) and the relation $[Q, {\bf b}(\theta)]=\ud(\theta)$ was also used. 
In order to prove the second relation of \rf{leib} one just need to use Jacobi identity for VOA. \footnote{We will also prove the generalized form of the Leibniz identity for 
$[,]$ in the next section using geometric arguments.  }
\hfill $\blacksquare$
\section{Homotopy Leibniz algebra for TCFT}
\noindent {\bf 3.0. Notation.} In this section, we will be dealing with full topological CFT (TCFT). By that we mean that we are working with the space $V\otimes \b V$,  where $V$ and $\bar V$ are TVOAs with essential elements $Q, J, N, L$ and $\bar{Q}, \bar{b}, \bar{N}, \bar{L}$. Whenever we write the expression $A(z, \b z)$ for the vector 
$A=\sum^n_{i=1}a_i\otimes {\bar a}_i$ so that $a_i\in V$,  ${\bar a}_i\in \bar{V}$, it means the following:  $A(z, \bar z)\equiv 
\sum^n_{i=1}a_i(z)\otimes {\bar a}_i(\bar z)$. Also we will use a shorthand notation: instead of $A( z, \bar z)$ we will write just $A(z)$, assuming dependence on 
$\bar{z}$ variable. 

 Most of the results below are valid for more generic full TCFTs. However, for simplicity we restrict ourselves to the case when the space of states is the tensor product of two VOAs.\\

\noindent {\bf 3.1. Definition of nonlocal operators in full TCFT.} To generalize the results of the section 2 in the case of full TCFT, one has to introduce new objects, which did not exist in the case of VOA, namely nonlocal operators. 

Let us consider the following expression under the correlator: 
\begin{eqnarray}\label{nl}
\int_V \langle... A_1^{(2)}(w_1)....A^{(2)}_n(w_n) B^{(1)}_1(v_1)....B^{(1)}_n(v_m)C(z)...\rangle. 
\end{eqnarray}
In the expression above $A^{(2)}_i(w)=A_i(w)dw\wedge d\b w$, $B^{(1)}_i(v)=B_i(v)dv+B'_i(v)d\b v$, $V$ is a compact manifold with boundary in the configuration space 
$\mathbb{C}^{n+m}\backslash D$, where $D$ is a collection of hyperplanes $\{w_i=w_j \}$, $\{v_k=w_l \}$, $\{v_r=v_s \}$, $\{v_i=\xi_s\}$, $\{w_i=\xi_p\}$, $\{v_i=\xi_q\}$, 
$\{w_i=z\}$, $\{v_r=z\}$ and the dots stand for other vertex operators, which are positioned in the points $\xi_1, \dots, \xi_d$.

Let us assume that $\max_{i,k}(|w_i-z|,|v_k-z|)=\rho$. Then we say that the resulting integrated object inside the correlator is the $nonlocal$ $operator$ at $z$ of size $\rho$. It is possible to shift it to  any other point $z'$ by means of standard translation operators $L_{-1} and \bar{L}_{-1}$.
One can consider all possible compositions of such operators, i.e. in \rf{nl} operators  $A_i(w)$, $B_k(v)$, $B'_k(v)$ could be also nonlocal, thus making the resulting nonlocal operator of possibly bigger size than original size $\rho$.  

We can treat these nonlocal operators on the same level as the standard vertex operators of full TCFT. For example, one can consider the following product 
under the correlator:
\begin{eqnarray}\label{nlcor}
 \langle A_1(z_1) A_2(z_2)...A_n(z_n)\rangle.
\end{eqnarray}
Here $A_1, ..., A_n$ are nonlocal operators of sizes $\rho_1, ..., \rho_n$. As long as $|z_i-z_j|<\rho_i+\rho_j$ for all $i, j$,  the expression \rf{nlcor} is well defined. 

Now we introduce an analogue of the space $W$ of differential forms in the case of full CFT. We will consider objects from $W\otimes \bar W$, so that $W, \bar W$ stands for spaces of differential forms associated with VOAs $V, \bar V$ where the formal odd variable from $W$ and $\bar W$ are denoted by $\theta$ and $\bar \theta$ correspondingly. Then this new space can be extended by nonlocal operators. In other words we are dealing with objects of three types, i.e. 0-forms, 1-forms and 2-forms:
\begin{eqnarray}\label{012}
\phi^{(0)}(z)=A(z),\quad \phi^{(1)}(z)=B(z)\theta+B'(z)\bar \theta, \quad \phi^{(2)}(z)=C(z)\theta\bar\theta,
\end{eqnarray}
where $A(z)$, $B(z)$, $B'(z)$, $C(z)$ are some nonlocal operators. Let us define the space $\mathcal{W}$ as a space spanned by all such objects. 

In order to define the analogues of operations $\circ$, $[\cdot,\cdot]$ on $\mathcal W$ we have to define integrals of 1-forms and 2-forms from \rf{012} over 1-dimensional or 2-dimensional compact manifolds  with or without boundary in $\mathbb{C}$, i.e. 
$\int_K\phi^{(1)}=\int_K  B(z)dz+B'(z)d\bar z$ and $\int_S\phi^{(2)}=\int_S C(z)dz\wedge d\bar z$, where $K$ and $S$ are respectively 1- and 2- dimensional submanifolds in $\mathbb{C}$.  
Let us define also the operator $\mathcal{D}$ acting on $\mathcal{W}$:
\begin{eqnarray}
\mathcal{D}=\mathcal{Q}+\ud, 
\end{eqnarray}
where $\mathcal{Q}=Q+\bar Q$ and $\ud=L_{-1}\theta+\bar{L}_{-1}\bar \theta$.  Another important object is the even operator ${\bf b}(z)=b_{-1}\theta+\bar{b}_{-1} \bar \theta$, such that 
$[\mathcal{Q}, {\bf b}]=\ud$.
\\

\noindent {\bf 3.2. Homotopy Leibniz identity.} Now we are ready to introduce an analogue of the operation $\circ$ in the case of $\mathcal{W}$. If $\phi , \chi \in \mathcal{W}$, then this new operation $\circ_{\epsilon}$
is defined as follows:
\begin{eqnarray}\label{circlep}
&&\phi\circ_{\epsilon} \chi(z,\theta)=\nonumber\\
&&\int_{C_{\epsilon,z}} \phi(w,\theta')\chi(z,\theta)-\int_{C_{\epsilon,z}}[\phi(w,\theta'),{\bf b}(\theta')]\chi(z,\theta),
\end{eqnarray}
where $\epsilon$ is the length of the radius of circle $C_{\epsilon,z}$ centered at $z$. we note that if only one term in the summand \rf{circlep} contributes to the expression: if $\phi$ is a 0-form, the first term is 0, if $\phi$ is a 1-form, the second term is zero, if $\phi$ is a 2-form this operation is equal to 0.
We also point out that $\circ_{\epsilon}$ is nonsingular only for those nonlocal $\phi$, $\chi$, the sizes of which ($\rho_{\phi}$, $\rho_{\chi}$) satisfy following relation: $\rho_{\phi}+\rho_{\chi}<\epsilon$.

Since we defined the counterpart of $\circ$ on $\mathcal{W}$, we are ready to define an analogue of the bilenar operation $[,]$ from Section 2:
\begin{eqnarray}
&&[\chi,\psi]_{\epsilon}(z,\theta)=\nonumber\\
&&\chi\circ_{\epsilon} \psi(z,\theta)-[\chi,{\bf b}(\theta)]\circ_{\epsilon} \psi(z,\theta)+\frac{1}{2}[[\chi,{\bf b}(\theta)], {\bb}(\theta)]\circ_{\epsilon} \psi(z,\theta).
\end{eqnarray}
One has to be careful with arguments here: for example if $\chi$ is a 1-form, then the second term is: 
\begin{eqnarray}
-\int_{C_{\epsilon,z}} [\chi(w,\theta'), \bb(\theta)]\psi(z,\theta)
\end{eqnarray}
Let us prove the following statement, which is an analogue of the first relation in Theorem 2.1.\\

\noindent {\bf Proposition 3.1.} {\it The operator $D$ and the operation $[\cdot ,\cdot]_{\epsilon}$ satisfy the relation:
\begin{eqnarray}\label{parleib}
D[\chi, \psi]_{\epsilon}=[D\chi, \psi]_{\epsilon}+(-1)^{|\chi|+1}[\chi, D\psi]_{\epsilon}+\nonumber\\
\ud f_{\chi, \epsilon}\psi+(-1)^{\chi}f^1_{\chi, \epsilon}\ud \psi,\nonumber\\
\end{eqnarray}
where  $f^1_{\chi, \epsilon}\psi=2[\chi,\psi]_{\epsilon}$.}\\

\noindent {\bf Proof.} 
 The relation \rf{parleib} is equivalent to
\begin{eqnarray}\label{leibbar}
\bar{D}[\chi, \psi]_{\epsilon}=[\bar{D}\chi, \psi]_{\epsilon}+(-1)^{|\chi|+1}[\chi, \bar{D}\psi]_{\epsilon},
\end{eqnarray}
where $\bar{D}=Q-\ud$. 
To simplify the notation, let us suppress the $\epsilon$-dependence of $[\cdot ,\cdot]$ and $\circ$:
\begin{eqnarray}
&&\bar {\mathcal{D}}[\chi,\psi]_{\epsilon}=(Q-\ud)[\chi, \psi]=\nonumber\\
&&Q\chi\circ \psi+(-1)^{|\chi|+1}\chi\circ Q\psi-\ud(\chi\circ\psi)+\nonumber\\
&&[\ud,\chi]\circ \psi +\ud([\chi, \bb]\circ \psi)+[\ud, [\bb, \chi]]\circ\psi- \nonumber\\
&&[Q\chi,{\bf b}]\circ \psi-(-1)^{|\chi|+1}[\chi,{\bf b}]\circ Q\psi+\nonumber\\
&& \frac{1}{2} [[Q\chi, \bb], \bb]\circ \psi +\frac{1}{2}(-1)^{|\chi|+1}[[\chi, \bb], \bb]\circ Q\psi= \nonumber \\
&&[\bar {\mathcal{D}}\chi, \psi]+
(-1)^{|\chi|+1}[\chi, \bar {\mathcal{D}}\psi].
\end{eqnarray}
Thus the proposition is proven. \hfill $\blacksquare$\\

\noindent{\bf Remark.} The appearance of extra term $\ud f_{\chi}\psi-(-1)^{\chi+1}f_{\chi}\ud \psi$ in \rf{parleib} 
may seem misleading because of clear Leibniz property for $\bar{D}$ \rf{leibbar}. As we will see later, such terms are unavoidable and will appear in all higher identities.\\

The next relation we would expect, following the ideas of Section 2, is the Leibniz property for the bracket for the parameter-dependent operation we described, however it appears 
that it holds only up to homotopy, which corresponds to parameter-dependent trilinear operation.
In order to define this trilinear operation properly we have to make some extra notation. 

Let  $D_{r,z}$ be a disk of radius $r$, centered at $z$ . Consider four positive real parameters $\rho_1>\rho_2$, $\epsilon<<\rho_{1,2}$, $\xi<<\rho_{1,2}$.  
We are interested in the following region in $\mathbb{C}^2$ (see Fig. 1): 
$V_3(z)=\{(v, w): v\in D_{\rho_1, z}\backslash (D_{\xi, w}\cup D_{\epsilon, z}), |z-w|=\rho_2\}$. 
\begin{figure}[hbt] 
\centering           
\includegraphics[width=0.5\textwidth]{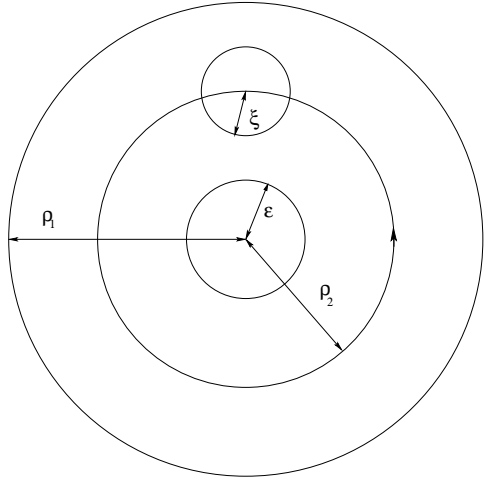}
\caption{Domain of integration $V_3$}                
\end{figure}   
Let us consider the following trilinear operation on $\mathcal{W}$:
\begin{eqnarray}
&&(\phi, \psi)_{\rho_1, \rho_2, \epsilon, \xi}\diamond \chi(z,\theta)=\frac{1}{2}\int_{V_3(z)}(\phi(v, \theta'')+[\phi(v, \theta''), \bb(\theta'')]+\nonumber\\
&&\frac{1}{2}[[\phi(v, \theta''), \bb(\theta'')], \bb(\theta''])(\phi(v, \theta')+[\phi(v, \theta'), \bb(\theta')])\chi(z,\theta).
\end{eqnarray} 
Then we can define the trilinear "bracket", which will be the desired homotopy for the Leibniz relation :
\begin{eqnarray}
&&[\alpha, \beta, \gamma]_{\rho_1, \rho_2, \epsilon, \xi}(z,\theta)=(\alpha, \beta)_{\rho_1, \rho_2, \epsilon, \xi}\diamond \gamma(z,\theta)-\nonumber\\
&&[(\alpha, \beta)_{\rho_1, \rho_2, \epsilon, \xi}, \bb(\theta)]\diamond \gamma(z,\theta)+
\frac{1}{2}[[(\alpha, \beta)_{\rho_1, \rho_2, \epsilon, \xi}, \bb(\theta)] ,\bb(\theta)] \diamond \gamma(z,\theta).
\end{eqnarray}

Then we have the following theorem.\\

\noindent{\bf Theorem 3.1.}{\it Operations $[\cdot, \cdot , \cdot]$ and $[, ]$ satisfy the following relation:
\begin{eqnarray}\label{jacobi}
&&D[\alpha,\beta,\gamma]_{\rho_1, \rho_2, \epsilon, \xi}-[D\alpha,\beta,\gamma]_{\rho_1, \rho_2, \epsilon, \xi}-\nonumber\\
&&(-1)^{|\alpha|}[\alpha,D\beta,\gamma]_{\rho_1, \rho_2, \epsilon, \xi}+(-1)^{|\alpha|+|\beta|}[\alpha,\beta,D\gamma]_{\rho_1, \rho_2, \epsilon, \xi}
=\nonumber\\
&&-(-1)^{|\alpha|}[[\alpha,\beta]_{\xi},\gamma]_{\rho_2}+(-1)^{|\alpha|}[\alpha,[\beta,\gamma]_{\rho_2}]_{\rho_1}-\nonumber\\
&&(-1)^{|\alpha|+(|\alpha|+1)(|\beta|+1)}[\beta,[\alpha,\gamma]_{\epsilon}]_{\rho_2}+\nonumber\\
&&\ud f^2_{(\alpha,\beta),V_3}\gamma+(-1)^{|\alpha|+|\beta|}f^2_{(\alpha,\beta),V_3}\ud\gamma,
\end{eqnarray}
where $f^2_{(\alpha,\beta), V_3}\gamma=2[\alpha,\beta,\gamma]_{\rho_1, \rho_2, \epsilon, \xi}$.}
\\

\noindent {\bf Proof.} Let us write explicitly term by term each of the terms from the RHS of \rf{jacobi}, suppressing the parameters for simplicity.
\begin{eqnarray}
&&-(-1)^{|\alpha|}[[\alpha,\beta],\gamma]=\\
&&-(-1)^{|\alpha|}[\alpha\circ \beta-[\alpha ,\bb]\circ \beta+\frac{1}{2}[[\alpha,\bb],\bb]\circ \beta ,\gamma]_{\rho_2}=\nonumber\\
&&-(-1)^{|\alpha|}((\alpha\circ\beta)\circ \gamma-[\alpha\circ \beta,\bb]\circ \gamma+\frac{1}{2}[[\alpha\circ\beta, \bb],\bb]\circ \gamma-\nonumber\\
&&([\alpha,\bb]\circ \beta)\circ \gamma+ [([\alpha,\bb]\circ\beta),\bb]\circ \gamma-\frac{1}{2}[[[\alpha,\bb]\circ\beta,\bb],\bb]\circ \gamma)=\nonumber\\
&&(-1)^{|\alpha|}((\alpha\circ\beta-[\alpha,\bb]\circ\beta)\circ \gamma-\nonumber\\
&&([[\alpha,\bb]\circ\beta-\alpha\circ\beta,\bb])\circ\gamma+\nonumber\\
&&\frac{1}{2}[[\alpha\circ\beta-[\alpha,\bb]\circ\beta, \bb],\bb]\circ\gamma),\nonumber
\end{eqnarray}
\begin{eqnarray}
&&(-1)^{|\alpha|}[\alpha,[\beta,\gamma]]=\\
&&(-1)^{|\alpha|}(\alpha\circ(\beta\circ\gamma)-[\alpha,\bb]\circ(\beta\circ\gamma)-\alpha\circ([\beta,\bb]\circ\gamma)+\nonumber\\
&&\frac{1}{2}([\alpha,\bb],\bb]\circ(\beta\circ\gamma)+[\alpha,\bb]\circ([\beta,\bb]\circ \gamma)+\frac{1}{2}\alpha\circ([[\beta, \bb], \bb]\circ \gamma),
\nonumber
\end{eqnarray}
\begin{eqnarray}
&&-(-1)^{|\alpha|+(|\alpha|+1)(|\beta|+1)}[\beta,[\alpha,\gamma]]=\\
&&-(-1)^{|\alpha|+(|\alpha|+1)(|\beta|+1)}(\beta\circ(\alpha\circ\gamma)-[\beta,\bb]\circ(\alpha\circ\gamma)-\beta\circ([\alpha,\bb]\circ\gamma)+\nonumber\\
&&\frac{1}{2}[\beta,\bb],\bb]\circ(\alpha\circ\gamma)+[\beta,\bb]\circ([\alpha,\bb]\circ \gamma)+\frac{1}{2}\beta\circ([[\alpha, \bb], \bb]\circ \gamma).\nonumber
\end{eqnarray}
Now one can compare this result with the expression from the LHS if one adds $f$-terms to it:
\begin{eqnarray}\label{st}
&&{\bar D}[\alpha,\beta,\gamma]_{\rho_1, \rho_2, \epsilon, \xi}-[D\alpha,\beta,\gamma]_{\rho_1, \rho_2, \epsilon, \xi}-\nonumber\\
&&(-1)^{|\alpha|}[\alpha,D\beta,\gamma]_{\rho_1, \rho_2, \epsilon, \xi}+(-1)^{|\alpha|+|\beta|}[\alpha,\beta,\bar{D}\gamma]_{\rho_1, \rho_2, \epsilon, \xi}
\end{eqnarray}
where $\bar{D}=Q-\ud$. 
If we apply Stokes theorem to \rf{st}, we obtain the proof.
\hfill $\blacksquare$\\

\noindent{\bf Remark.} One can notice that in the chiral case one can consider the algebraic structure of 0-forms and 1-forms based on operation $\circ$ and it appears to be again a Leibniz algebra with a differential $Q$. The Leibniz algebra of 0-forms reproduces the Lian-Zuckerman bracket \cite{lz}.
The analogues of these operations in the case of full CFT does not produce a Leibniz algebra or a Leibniz algebra up to homotopy as one can show by a direct calculation on the whole space $\mathcal{W}$. However, one can restrict 
$\circ$ to 1-forms only and, using the Stokes theorem, show that it satisfies Leibniz algebra relations up to homotopy for $D$ only. Unfortunately, $\circ$ does not satisfy higher homotopy Leibniz relations being restricted to 1-forms as we will see in the following.\\

\noindent {\bf 3.3. Higher Order Leibniz Identity in TCFT.}
In this subsection, we prove the higher order relation, involving bilinear operation $[\cdot,\cdot ]$ and trilinear one $[\cdot, \cdot, \cdot] $. It turns out that this generalized 
Leibniz identity again holds up to homotopy which is given by a quadrilinear operation. 
In order to define this quadrilinear operation properly let us introduce a shorthand notation. Let $\phi\in \mathcal{W}$, then we define  ${\bf \phi}$ as follows:
\begin{eqnarray}
{\hat \phi}(z,\theta)=\phi(z, \theta)+[\phi(z, \theta), \bb(\theta)]+\frac{1}{2}[[\phi(z, \theta), \bb(\theta)],  \bb(\theta)].
\end{eqnarray}
Let us consider the following region in $\mathbb{C}^3$ (see Fig. 2) centered at z: $V_4(z)=\{(u,v, w): u\in D_{z, \rho_1}\backslash(D_{\alpha, v}\cup D_{\xi_1,w}\cup D_{\epsilon_1,z});  
v\in D_{z, \rho_2}\backslash (D_{\xi_2,w}\cup D_{\epsilon_2,z}); |z-w|=\rho_3\}$. The relations between real positive parameters $\rho_i, \epsilon_k, \xi_l. \alpha$ are as follows:$\rho_i>>\alpha, \xi_k, \epsilon_l$ for all $i,k, l$. Also, $\epsilon_2>>\epsilon_1, \alpha$,  $\xi_2>>\xi_1, \alpha$. 

\begin{figure}[hbt] 
\centering           
\includegraphics[width=0.7\textwidth]{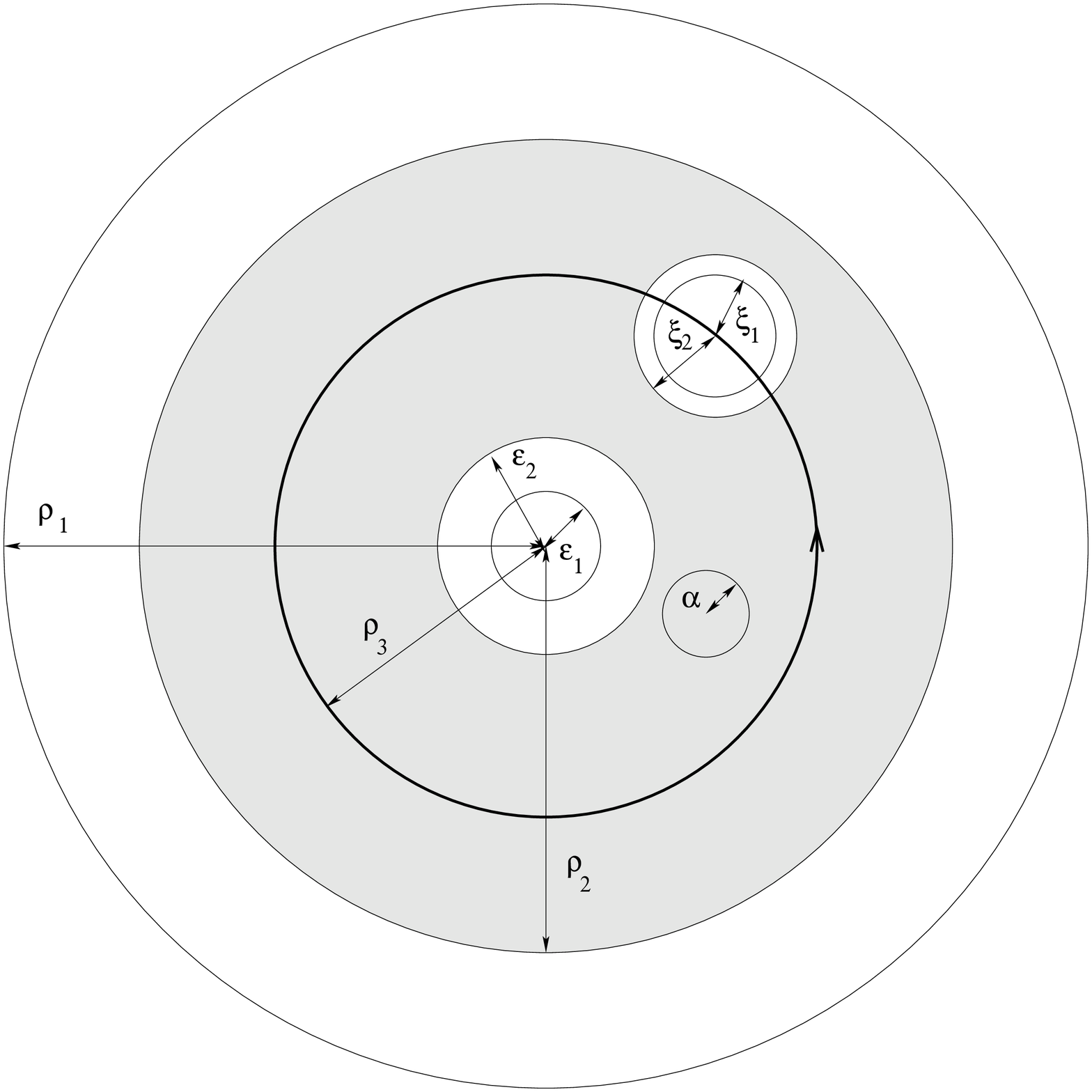}
\caption{Domain of integration $V_4$}                
\end{figure}

Using all this data we are ready to define the following quadrilinear operation: let $\phi, \psi, \chi, \sigma\in \mathcal{W}$, then 
\begin{eqnarray}
(\phi, \psi, \chi)_{V_4} \diamond \sigma(z,\theta)=\frac{1}{4}\int_{V_4(z)}\hat{\phi}(u,\theta''')\hat{\psi}(v,\theta'')\hat{\chi}(w,\theta')\sigma(z,\theta).
\end{eqnarray}
We define the operation $[\cdot, \cdot, \cdot, \cdot]$ which will be a homotopy for generalized Leibniz identity as follows:
\begin{eqnarray}
&&[\phi, \psi, \chi, \sigma]_{V_4}(z,\theta)=(\phi, \psi, \chi)_{V_4} \diamond \sigma(z,\theta)-\nonumber\\
&&[(\phi, \psi, \chi)_{V_4}, \bb(\theta)]\diamond \sigma(z,\theta)+\frac{1}{2}[[(\phi, \psi, \chi)_{V_4}, \bb(\theta)],\bb(\theta)] ] \diamond \sigma(z,\theta)
\end{eqnarray}

Now we are ready to formulate a theorem.\\

\noindent{\bf Theorem 3.2.} {\it Let $\phi, \psi, \chi, \sigma\in \mathcal{W}$. Then the operations $[\cdot, \cdot, \cdot, \cdot]$, $[\cdot, \cdot, \cdot]$, $[\cdot, \cdot]$ satisfy the following relations:
\begin{eqnarray}\label{genjacobi}
&&\mathcal{D}[\phi, \psi, \chi, \sigma]_{V_4}-[\mathcal{D}\phi, \psi, \chi, \sigma]_{V_4}-(-1)^{|\phi|}[\phi,\mathcal{D} \psi, \chi, \sigma]_{V_4}\label{bound1}\\
&&-(-1)^{|\phi|+|\psi|}[\phi, \psi,\mathcal{D} \chi, \sigma]_{V_4}+(-1)^{|\phi|+|\psi|+|\chi|}[\phi,\psi,\chi,\mathcal{D} \sigma]_{V_4}=\nonumber\\
&&-(-1)^{|\phi|}((\phi,\psi)_{\alpha},\chi,\sigma)_{\rho_2,\rho_3,\epsilon_2,\xi_2}+(-1)^{|\phi|}(\phi,(\psi,\chi,\sigma)_{\rho_2,\rho_3,\epsilon_2,\xi_2})_{\rho_1}-\label{first4}\\
&&(-1)^{|\phi|+(|\phi|+1)(|\psi|+|\chi|+1)}(\psi,\chi, (\phi,\sigma)_{\epsilon_1})_{\rho_2,\rho_3, \epsilon_2,\xi_2}-\nonumber\\
&&(-1)^{|\phi|+(|\phi|+1)|\psi|}(\psi, (\phi,\chi)_{\xi_1},\sigma)_{\rho_2,\rho_3, \epsilon_2,\xi_2}+\nonumber\\
&&-(-1)^{|\phi|+|\psi|+(|\psi|+1)|\phi|}(\psi, (\phi, \chi, \sigma)_{\lambda_1,\rho_3, \xi_1, \epsilon_1})_{\rho_2}+\label{l1}\\
&&-(-1)^{|\phi|+|\psi|} (\phi, \psi, (\chi, \sigma)_{\rho_3})_{\rho_1,\rho_2, \lambda_1, \alpha}-\nonumber\\
&&(-1)^{|\phi|+|\psi|+(|\chi|+1)(|\phi|+|\psi|+1)}(\chi, (\phi, \psi, \sigma)_{\lambda_2, \epsilon_2, \alpha, \epsilon_1})_{\rho_2}-\label{l2}\\
&&(-1)^{|\phi|+|\psi|+(|\chi|+1)(|\psi|+1)}(\phi, \chi, (\psi, \sigma)_{\epsilon_2})_{\rho_1,\rho_3, \lambda, \alpha}-\nonumber\\
&&(-1)^{|\phi|+|\psi|}((\phi,\psi,\chi)_{\lambda_3,\xi_2, \alpha,\xi_1},\sigma)_{\rho_2}- (-1)^{|\phi|+|\psi|}(\phi, (\psi,\chi)_{\xi_2},\sigma)_{\rho_1,\rho_3,\lambda_3, \epsilon_1}\label{l3}+\nonumber\\
&&\ud f^3_{(\phi,\psi,\chi), V_4}\sigma+(-1)^{|\phi|+|\psi|+|\chi|}f^3_{(\phi,\psi,\chi), V_4}\ud \sigma,\label{bound2}
\end{eqnarray}
where $f^3_{(\phi,\psi,\chi), V_4}\sigma=2[\phi, \psi, \chi, \sigma]_{V_4}$. 
The conditions on the real parameters $\lambda_i$ are as follows: $\rho_3+\epsilon_1<\lambda_1\rho_2+\alpha$, $\rho_3+\epsilon_1>\lambda_2>\epsilon_2+\alpha$, $\rho_3-\epsilon_1>\lambda_3>\epsilon_2+\alpha$.}\\

\noindent {\bf Proof.} To prove it one has to use Stokes theorem agai, combining terms \rf{bound1}, \rf{bound2}, now for the region $V_4$. Let us describe in detail how all the terms \rf{first4}, \rf{l1}, \rf{l2}, \rf{l3} arise from the boundary of $V_4$. The first four terms \rf{first4} correspond to the four boundary circles of the outer disc  $D_{z, \rho_1}\backslash(D_{\alpha, v}\cup D_{\xi_1,w}\cup D_{\epsilon_1,z})$. However, each of pairs \rf{l1}, \rf{l2}, \rf{l3} comes from each of the boundary components of $D_{z, \rho_2}\backslash (D_{\xi_2,w}\cup D_{\epsilon_2,z})$. Let us write down explicitly how each of these six terms appear.
First of all, let us consider the outer boundary of $D_{z, \rho_2}\backslash (D_{\xi_2,w}\cup D_{\epsilon_2,z})$. The corresponding boundary piece of $V_4$ is given on Fig. 3. In order to represent it via $V_3$ we divide  
$D_{z, \rho_2}$ into two pieces by means of circle of the radius $\lambda_1$ such that $\rho_3+\epsilon_1<\lambda_1\rho_2+\alpha$.
\begin{figure}[hbt] 
\centering           
\includegraphics[width=0.6\textwidth]{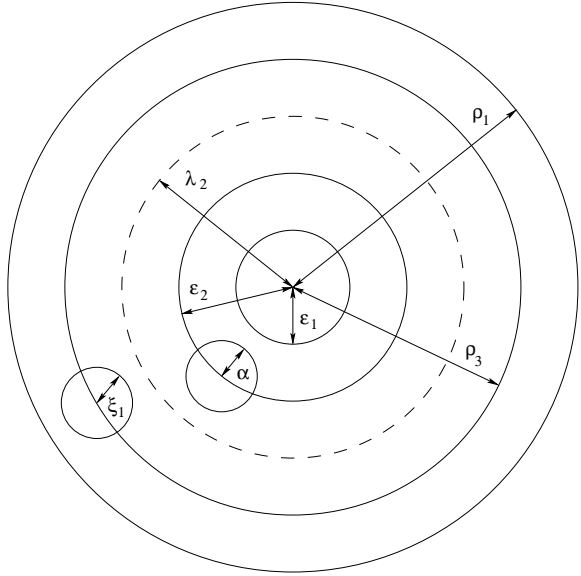}
\caption{Boundary component of $V_4$}                
\end{figure}
Similarly we decompose boundary components corresponding to inner boundaries, i.e. circle of the radius $\epsilon_2$
with center at $z$ (see Fig. 4)
\begin{figure}[hbt] 
\centering           
\includegraphics[width=0.6\textwidth]{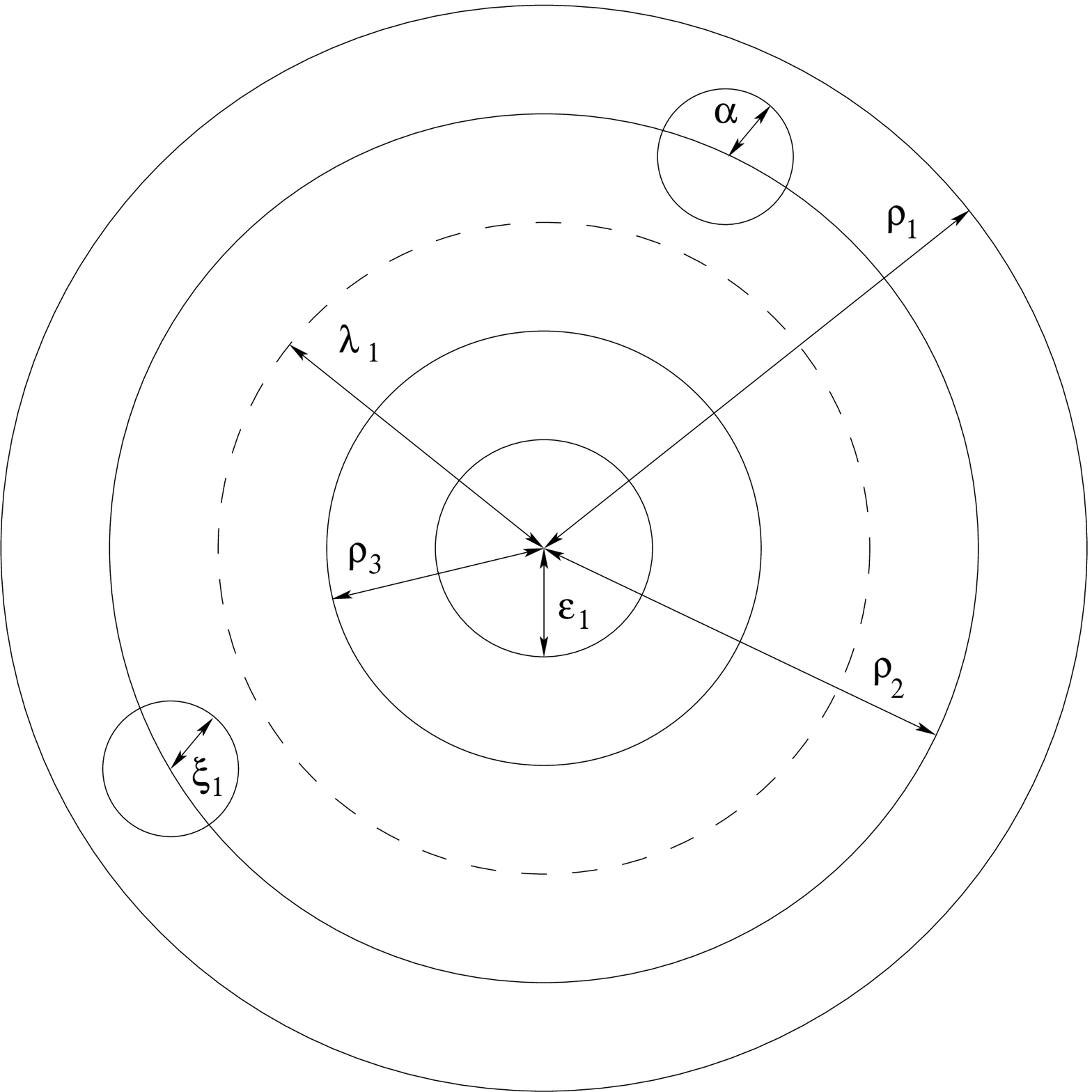}
\caption{Boundary component of $V_4$}                
\end{figure}
and the circle of the radius $\xi_2$ with the center at $w$, lying on the circle of radius $\rho_3$ (see Fig. 5).
\begin{figure}[hbt] 
\centering           
\includegraphics[width=0.6\textwidth]{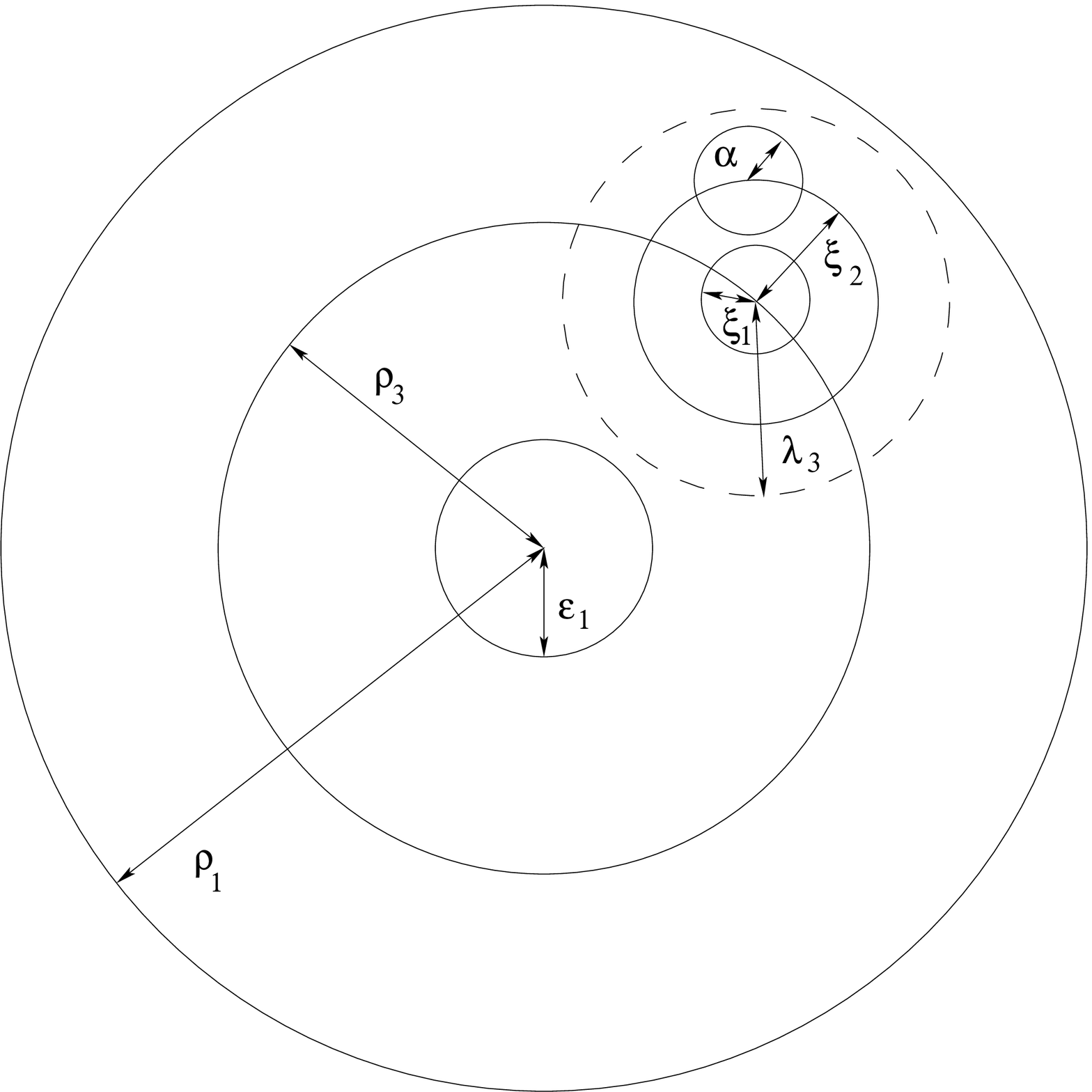}
\caption{Boundary component of $V_4$}                
\end{figure}
Integration over these six boundary components gives the necessary terms  \rf{l1}, \rf{l2}, \rf{l3}.

\hfill$\blacksquare$\\

\section{Conjectures}

\noindent {\bf 4.1. Loday's homotopy Leibniz algebras and higher operations for TCFT.} The results of the previous section suggest that there could be more of higher order $n$-linear operations on $\mathcal{W}$, such that they satisfy higher Leibniz-like identities. The relations we already obtained \rf{parleib}, \rf{jacobi}, \rf{genjacobi} up to the $f$-terms and dependence of parameters are the first four relations for the Loday's $Leibniz_{\infty}$.  

Recall that a $Leibniz_{\infty}$ algebra \cite{loday}, \cite{xu} is a $\mathbb{Z}$-graded vector space with the set of multilinear maps $(\mu_k)^{\infty}_{k=1}$,  $\mu_k:\otimes^k V\to V$ of degree $1$ satisfying the following identities:
\begin{eqnarray}\label{lod}
&&\sum_{1\le j,k\le n}\sum_{\mathcal{S}^{j-1}_{k-j}}\epsilon(\sigma, \phi_1, ..., \phi_{k-1})(-1)^{|\phi_{\sigma(1)}|+
|\phi_{\sigma(2)}|+\dots+|\phi_{\sigma(k-j)}|}\nonumber\\
&&\mu_{n-j+1}(\phi_{\sigma(1)}, \dots, \phi_{\sigma(k-j)}, \mu_{j}(\phi_{k+1-j}, \dots, \phi_{k}), \phi_{k+1},..., \phi_n)=0,
\end{eqnarray}
where $n\in \mathbb{N}$ and homogeneous $\phi_i\in V$. Here $\mathcal{S}^{j-1}_{k-j}$ denotes the $(k-j, j-1)$ shuffled permutations and $\epsilon(\sigma, \phi_1, ..., \phi_{k-1})$ is the Koszul sign of the permutation of 
$\phi_1,\dots ,  \phi_k$.

Identities between $\mu_1$, $\mu_2$, $\mu_3$, $\mu_4$ are the relations between $\mathcal{D}$, $[\cdot, \cdot]$, $[\cdot,\cdot, \cdot]$, $[\cdot,\cdot, \cdot, \cdot]$ modulo $f$-terms, dependence of parameters and appropriate shift in grading.

Therefore, natural conjecture is that the higher Leibniz relations in our case will be like in \rf{lod} modulo the $\ud$-homotopy term, i.e. if we denote LHS of \rf{lod} by $(Leibniz_{\infty}-{\rm relations})(\phi_1,\cdot, \phi_n)$ we expect the following
\begin{eqnarray}\label{leibc}
&&(Leibniz_{\infty}-{\rm relations})(\phi_1,\dots, \phi_n)=\nonumber\\
&&\ud f^{n-1}_{(\phi_1,..., \phi_{n-1})}\phi_n+(-1)^{|\phi_1|+...+|\phi_{n-1}|+1}f^{n-1}_{(\phi_1,..., \phi_{n-1})}\ud \phi_n,
\end{eqnarray}
where $f^{n-1}_{(\phi_1,..., \phi_{n-1})}\phi_n=2\mu_n(\phi_1,..., \phi_{n})$.

Natural interesting problem is the construction of the explicit form of the generalized $\mu_n$. It appears that the compact manifolds with boundary $V_3$ and $V_4$ considered in subsections 3.2 and 3.3. are topologically equivalent to compactified moduli spaces $\bar{M}(3)$ and $\bar{M}(4)$. Here $\bar{M}(n)$ is the compactification of moduli space  $\mathbb{C}^{n}\backslash\{\rm diag\}/G$ ($G$ is the group consisting of translations $\mathbb{R}^2$ and dilations $\mathbb{R}_+$).

Therefore, the proper way to construct the multilinear operations is as follows. First, as in subsections 3.2 and 3.3 we define the following n-linear operations:
\begin{eqnarray}
(\hat{\phi}_1, ..., \hat{\phi}_{n-1})_{V_n}\diamond\phi_n(z,\theta)=\int_{\bar{M}_z(n)}\hat{\phi}_1\dots \hat{\phi}_{n-1}\phi_{n}(z,\theta),
\end{eqnarray}  
where compactified moduli space $\bar{M}_z(n)$ is centered around point $z$ and depends on the set of real parameters, and $V_2\equiv M(2)$ is just a circle around $z$. Then, in order to define multilinear bracket $[\cdot, \dots, \cdot]_n$ and therefore its degree 1 counterpart $\mu_n$, one has to add the action of $b(\theta)$ on the first $n-1$ variables:
\begin{eqnarray}
&&[\phi_1, ..., \phi_{n-1},\phi_n]_n(z,\theta)=(\hat{\phi}_1, ..., \hat{\phi}_{n-1})\diamond\phi_n(z,\theta)-\nonumber\\
&&[(\hat{\phi}_1, ..., \hat{\phi}_{n-1}),\bb(\theta)]\diamond\phi_n(z,\theta)+\frac{1}{2}[[(\hat{\phi}_1, ..., \hat{\phi}_{n-1}),\bb(\theta)],\bb(\theta)]\diamond\phi_n(z,\theta). 
\end{eqnarray}  

The relations \rf{leibc} should be proved by means of the cell decompositions of $\bar{M}(n)$.\\

\noindent{\bf 4.2. Generalized Maurer-Cartan equation, $\beta$-function and physical interpretation.} In the subsection 4.1, we conjectured  the existence of the parameter-dependent homotopy Leibniz algebras on the space $\mathcal{W}$. In general, for standard $Leibniz_{\infty}$ algebras one can write a generalized Maurer-Cartan equation associated to it. Explicitly this equation  is as follows:
\begin{eqnarray}
D\phi+\sum^{\infty}_{n=2}[\phi,\dots,\phi]_n=0,
\end{eqnarray} 
where $\phi$ is an element of degree 2. 
As well as the standard Maurer-Cartan equation, this equation has symmetries generated by the elements of degree $1$. We have some obstacles for construction of Maurer-Cartan equation in the case of the algebras we want to construct. 
One obstacle is that our 
operations $[\cdot,..., \cdot]_n$ are parameter-dependent and another is that $Leibniz_{\infty}$ relations hold up to $f$-terms. In order to get rid of the $f$-terms we will change the RHS of the Maurer-Cartan equation so that it is satisfied up to homotopy with respect to the de Rham operator $\ud$. As for parameter dependence, we conjecture that if we consider arguments of the multilinear operations $[\cdot,..., \cdot]_n$ to be local, one can expand the operations in terms of the parameters of $\mathcal{M}(n)$. The corresponding 0-modes in this expansion will satisfy the required homotopy Leibniz algebra relations without parameter-dependence. 

It was shown in \cite{lmz}, \cite{zeit2}, \cite{zeit3}, that this parameter expansion is possible for bilinear operation, and that the Maurer-Cartan equation expanded up to the second order for some CFTs, reproduces 
what physicists know as $\beta$-functions for sigma models, if $\mathcal{Q}$ is the semi-infinite cohomology (BRST) operator. 

It is natural to interpret 0-forms, 1-forms and 2-forms as correspondingly local obeservables, currents and perturbations. The presence of $f$-terms in algebraic relations corresponds to the ambiguity in defining currents and perturbations in peturbative physical theories: currents and perturbations are usually defined up to exact term with respect to de Rham differential.

We claim that the Maurer-Cartan equation for the $Leibniz_{\infty}$ algebra (up to homotopy) we discovered in this paper gives the proper description of the $\beta$-function in perturbed CFT.

We plan to return to these and many other related questions, including higher dimensional generalization in the forthcoming articles.

\end{document}